\documentclass[preprint,12pt]{elsarticle}

\newtheorem{theorem}{Theorem}[section]
\newtheorem{lemma}{Lemma}[section]
\newtheorem{corollary}{Corollary}[section]
\newtheorem{proposition}{Proposition}[section]

\newcommand{\ignore}[1]{}{}

\usepackage{amssymb}
\usepackage{amsmath}

\usepackage{color}
\usepackage{soul}
\usepackage[dvipsnames]{xcolor}

\usepackage[colorlinks=true, linkcolor=blue, citecolor=blue]{hyperref}

\def\1{{{\mbox{${\rm{1\negthinspace\negthinspace I}}$}}}}

\renewcommand{\theequation}{\arabic{equation}}

\newcommand\beq{\begin{equation}}
\newcommand\eeq{\end{equation}}
\renewcommand{\theequation}{\thesection.\arabic{equation}}

\usepackage{ifthen}
\usepackage{xkeyval}
\usepackage{todonotes}
\begin{document}

\begin{frontmatter}

\title{Uniform  Cram\'{e}r moderate deviations  and Berry-Esseen bounds for a supercritical branching process in a random environment}
\author[cor1]{Xiequan Fan$^*$}
\author[cor2]{Haijuan Hu}
\author[cor3]{Quansheng Liu}
 \cortext[cor1]{\noindent Corresponding author. \\
\mbox{\ \ \ \ }\textit{E-mail}:  fanxiequan@hotmail.com (X. Fan)\\
\mbox{\ \ \ \ \ \ \ \ \ \ \  \ \ \ \ } . }
\address[cor1]{Center for Applied Mathematics,
Tianjin University, Tianjin 300072, China }
\address[cor2]{School of Mathematics and Statistics, Northeastern University at Qinhuangdao, Qinhuangdao, China}
\address[cor3]{Universit\'{e} de Bretagne-Sud, LMBA, UMR CNRS 6205, Campus de Tohannic,  56017 Vannes, France}

\begin{abstract}
Let $\{Z_n, n\geq 0\}$ be a supercritical branching process in an independent and identically distributed random environment.
We prove Cram\'{e}r  moderate deviations and Berry-Esseen bounds for $\ln (Z_{n+n_0}/Z_{n_0})$ 
uniformly in $n_0 \in \mathbb{N}$, which extend
the corresponding results
by  Grama et al. 
 (Stochastic Process.\ Appl. 2017) 
established for $n_0=0$. 
The extension is 
interesting in theory,  and is motivated by applications.
A new method is developed for the proofs; some conditions of Grama et al. (2017) are relaxed in our present setting.
 An example of  application is given in  constructing  confidence intervals to estimate  the criticality parameter in terms of  $\ln (Z_{n+n_0}/Z_{n_0})$ and $n$.
\end{abstract}

\begin{keyword} Branching processes; Random environment;  Cram\'{e}r  moderate deviations; Berry-Esseen bounds
\vspace{0.3cm}
\MSC primary 60J80; 60K37; secondary  60F10
\end{keyword}

\end{frontmatter}




\section{Introduction}
\setcounter{equation}{0}
As an important generalization  of the Galton-Watson process, the branching process in a random environment (BPRE)
was first  introduced by Smith and Wilkinson \cite{SW69} to modelize the growth of a population submitted to an independent and identically distributed (iid) random environment. Basic results for a BPRE can be found in   Athreya and Karlin \cite{AK71a,AK71b} who considered  the stationary and ergodic environment case.

 A BPRE can be described as follows.  Let $\xi=(\xi_0, \xi_1, ...) $ be a sequence of independent and identically distributed (iid)  random variables,
 where  
  $\xi_n$ stands for the random
environment at time $n.$  Each realization of $\xi_n$ corresponds to a probability law $ p(\xi_n) = \{ p_i(\xi_n): i \in  \mathbb{N}\}$  on
 $\mathbb{N}=\{0, 1,...\}$  ($p_i(\xi_n) \geq 0  $ and $ \sum_{i=0}^\infty p_i(\xi_n)=1$).
 A branching process $\{Z_n, n\geq 0 \} $ in the random environment $\xi$ can be defined as follows:
$$Z_0=1,\ \ \ \ \ \ \ \  Z_{n+1}= \sum_{i=1}^{Z_n} X_{n,i} \   \ \ \ \textrm{for }\   n \geq 0,$$
where $X_{n,i}$ is the number of offspring  of the $i$-th individual in generation $n.$
Conditioned
on the environment $\xi$ 
the random variables $X_{n,i} \,  (n\geq 0, i\geq 1)$ are independent, and  each $X_{n,i}$ has the   same law  $p(\xi_n)$.
Denote by $\mathbb{P}_\xi$  the  probability when the environment
$\xi$ is given, $\tau$  the law of the environment $\xi$, and $$\mathbb{P}(dx, d\xi )=\mathbb{P}_\xi(dx)\tau(d\xi)$$
the total law of the process; $\mathbb{P}_\xi$ can be considered as the conditional law of $\mathbb{P}$ given the environment $\xi$.
The  conditional probability $\mathbb{P}_\xi$ is called the quenched law,  while the total probability $\mathbb{P}$ is
 called annealed law. In the sequel  $\mathbb{E}_\xi$ and $\mathbb{E}$ denote respectively the quenched and annealed expectations.
Set  for $ n\geq 0$,
$$m_n   = \sum_{k=0}^{\infty} k \, p_k(\xi_n)  \quad  \mbox{ and  }  \quad  \Pi_{n}= \prod_{i=0}^{n-1} m_{i},  $$
with the convention that $ \Pi_{0}= 1$. Then
$ m_n =\mathbb{E}_\xi X_{n, i}$ for each $i\geq 1$ and  $ \Pi_{n}=  \mathbb{E}_\xi Z_{n}$. Let
$$X= \log m_0,\ \ \ \mu =\mathbb{E}X.  $$ 
 The process $\{Z_n, n\geq 0 \}$  is called supercritical, critical or subcritical according to $\mu >0$, $\mu =0$ or $\mu <0$, respectively. We call $\mu$ the criticality parameter.

 Limit theorems for BPRE have attracted a lot of attentions.
 See for example Vatutin \cite{V10}, Afanasyev et al.\ \cite{ABV14}, Vatutin and Zheng \cite{VZ12}  and  Bansaye  and  Vatutin \cite{BV17}
on the survival probability and conditional limit theorems for  subcritical BPRE.
 For supercritical  BPRE,   a number of researches have studied  moderate and large deviations; see, for instance, \mbox{Kozlo} \cite{K06},
Bansaye  and Berestycki \cite{VB09},  B\"{o}inghoff  and Kersting \cite{BK10}, Bansaye and  B\"{o}inghoff \cite{BB11},
  Huang and Liu \cite{HL12}, 
 Nakashima \cite{N13},   B\"{o}inghoff \cite{B14}, and Grama, Liu and Miqueu\ \cite{GLE17}.

In this paper, we are interested in  Cram\'{e}r moderate deviations and Berry-Esseen bounds for a supercritical  BPRE.
For simplicity we assume that
\begin{equation} \label{cond-0}
p_0(\xi_0)=0 \quad \mathbb{P}\textrm{-a.s.} \quad \mbox{ and  } \quad \sigma^2=\mathbb{E}(X-\mu)^2 \in (0, \infty),
\end{equation}
which imply  that the process is supercritical and $Z_n \rightarrow \infty$  a.s. 
Under the additional conditions:
$ \mathbb{E}  \frac{Z_1^{p}}{m_0}   < \infty$ for a constant $ p>1$ and    $\mathbb{E}e^{\lambda_0 X}   < \infty$
for a constant $\lambda_0>0$,
  Grama et al.\ \cite{GLE17} have established the Cram\'{e}r   moderate deviation expansion, which implies in particular that
for $0 \leq x = o(\sqrt{n} )$ as $n\rightarrow \infty$,
\begin{equation}\label{cramer}
\Bigg| \ln \frac{\mathbb{P}\big( \frac{ \ln  Z_{n}  - n \mu \ }{  \sigma \sqrt{n}} \geq x  \big)}{1-\Phi(x)} \Bigg|  \leq C    \frac{  1+x^3    }{  \sqrt{n}  },
\end{equation}
where  throughout the paper the symbol  $C$,  probably  supplied  with some indices,
denotes  a positive  constant whose value may differ from line to line.
Inequality (\ref{cramer}) is interesting due to the fact that it implies a moderate  deviation principle  (MDP)  and the following result about the equivalence to the normal tail:
\begin{equation}\label{cramedr}
\frac{\mathbb{P}\big( \frac{ \ln  Z_{n}  - n \mu \ }{  \sigma \sqrt{n}}   \geq x  \big)}{1-\Phi(x)} =1+o(1),
\end{equation}
for  $x \in [0, o(n^{1/6}))$, as $n\rightarrow \infty$.
Assuming $ \mathbb{E}  \big(\frac{Z_1}{m_0} \big )^{p} < \infty$ for a constant $ p>1 $ and   $\mathbb{E}X^{2+\rho}    < \infty$
  for a constant $\rho \in (0, 1)$,  Grama et al.\ \cite{GLE17}   have also obtained the following Berry-Esseen bound  for $\ln Z_n$:
\begin{equation} \label{bercrser}
\sup_{x\in \mathbb{R}} \Big|\mathbb{P}\Big( \frac{ \ln  Z_{n}  - n \mu \ }{  \sigma \sqrt{n}} \leq  x  \Big)  -  \Phi(x)   \Big|
\leq  \frac{  C }{   n^{\rho/2}   }.
\end{equation}

The results  (\ref{cramer}), (\ref{cramedr}) and (\ref{bercrser}) are interesting both in theory and in applications. For example, they can be applied to obtain confidence intervals  to estimate the criticality parameter $\mu$ in terms of the observation $Z_n$ and the present time $n$,  or to estimate the population size $Z_n$ in terms of $\mu $ and $n$.
In the real-world applications, 
 it may happen that we know a historical data $Z_{n_0}$ for some $n_0>0$, the current population size $Z_{n_0+n}$,  as well as the increment $n$ of generation numbers,
 but do not know the generation number $n_0 + n$.  In such a case (\ref{cramer}), (\ref{cramedr}) and (\ref{bercrser}) are no longer applicable to obtain confidence intervals
 to estimate $\mu$ in terms of $Z_{n_0}, Z_{n_0+n}$  and $n$, while $n_0>0$.
 The same problem exists while we want to construct confidence intervals to preview  $Z_{n_0+n}$  in terms of $Z_{n_0}, \mu$  and $n$.
 Motivated by these problems, we will
  extend (\ref{cramer}), (\ref{cramedr})   and (\ref{bercrser}),  with $ \ln  Z_{n}$   replaced by $ \ln  \frac{Z_{n_0+n} }{Z_{n_0}}$, uniformly in  $n_0 \in \mathbb N$ (so that in applications $n_0$ can be  taken as an function of $n$).   This is the main objective of the present paper.

The main results are presented in Section \ref{sec2}. Let us introduce them briefly.
Denote by $x^+=\max\{x, 0\}$ the positive   part  of $x$.
  In Theorem \ref{th2.2},
assuming  $\mathbb{E} \frac{Z_1}{m_0} \ln^+ Z_1  < \infty $ and
 $\mathbb{E}e^{\lambda_0 X}   < \infty$ for a constant $\lambda_0>0$,
 we prove that uniformly in $n_0 \in \mathbb N$,  for     $0 \leq x  = o(\sqrt{n } )$,  as $ n \rightarrow \infty$,
\begin{equation} \label{th02s}
 \Bigg| \ln \frac{\mathbb{P}\big( \frac{ \ln  \frac{Z_{n_0+n}}{Z_{n_0}}  - n \mu \ }{  \sigma \sqrt{n}} \geq x  \big)}{1-\Phi(x)} \Bigg|
  \leq C      (1+x^3 ) \frac{ 1+  \mathbf{1}_{[0, \sqrt{  \ln n }) } (x) \ln n  }{  \sqrt{n}  }.
\end{equation}
When $n_0=0$, inequality (\ref{th02s}) reduces nearly to (\ref{cramer}), with    $\ln n$ as an additional factor.
Notice that here we do not need the additional condition  that $ \mathbb{E}  \frac{Z_1^{p}}{m_0}   < \infty$ for some $ p>1$ assumed in  
 \cite{GLE17} for  (\ref{cramer}) to hold.
 As a consequence, we obtain a uniform MDP for  
$ \ln  \frac{Z_{n_0+n}}{Z_{n_0}}$, see
Corollary \ref{corollary02}.  From  \eqref{th02s}, we also obtain
the following equivalence to the normal tail: uniformly in $n_0 \in \mathbb N$, for  $x \in [0, o(n^{1/6}))$, as $n\rightarrow \infty$,
\begin{equation}\label{cramedr002}
\frac{\mathbb{P}\big( \frac{ \ln    \frac{Z_{n_0+n}}{Z_{n_0}}    - n \mu \ }{  \sigma \sqrt{n}}   \geq x  \big)}{1-\Phi(x)} =1+o(1).
\end{equation}
 When the exponential moment condition $\mathbb{E}e^{\lambda_0 X}   < \infty$ is relaxed to the sub-exponential moment condition that
 $\mathbb{E}\exp\{  \lambda_0  X^{\frac{4 \gamma}{1-2\gamma}} \}<\infty$ for some  $\gamma \in (0, \frac16]$,
we prove that \eqref{cramedr002} still holds  for $x \in [0, o(n^{\gamma}))$; see  Theorem \ref{th03} for a result of type Linnik \cite{L61}.
Using  \eqref{cramedr002}, we can prove, under  the exponential moment condition,
the following uniform Berry-Esseen bound: 
uniformly in $n_0 \in \mathbb N$,
 \begin{equation} \label{bercrser002}
\sup_{x\in \mathbb{R}} \Big|\mathbb{P}\Big( \frac{ \ln   \frac{Z_{n_0+n}}{Z_{n_0}}  - n \mu \ }{  \sigma \sqrt{n}} \leq  x  \Big)  -  \Phi(x)   \Big|
\leq  \frac{  C \ln n  }{   \sqrt n   }.
\end{equation}
  Compared to the best rate  $\frac{ C}{  \sqrt{n}  }$ of the Berry-Esseen bound  for random walks,   here the  factor $\ln n$ is added.
  We believe that this  factor $\ln n$ can be removed from \eqref{bercrser002},
   just as in the case  $n_0 = 0$ considered in Grama et al.\ \cite{GLE17}.
In fact  for $n_0=0$, the more general Berry-Esseen bound $ \frac{  C }{   n^{\rho/2}}$ was established   in  \cite{GLE17}
under  the moment condition $\mathbb{E}X^{2+\rho}    < \infty$ with
$\rho \in (0, 1] $. In this paper,  we prove that
   if
  $\mathbb{E}X^{2+\rho}    < \infty$ for some $ \rho \in (0, \frac{\sqrt{5} -1}{2}), $ then uniformly in $n_0 \in \mathbb N$,
\begin{equation} \label{main01}
\sup_{x\in \mathbb{R}} \Big|\mathbb{P}\big(  Z_{n_0,n} \leq  x  \big)  -  \Phi(x)   \Big|
\leq  \frac{  C }{   n^{\rho/2}   }.
\end{equation}
See Theorem \ref{th02}. Clearly, inequity (\ref{main01}) with $n_0=0$ reduces to  (\ref{bercrser}), which was obtained in \cite{GLE17} under the additional condition that
$ \mathbb{E}  \big(\frac{Z_1}{m_0} \big )^{p} < \infty$ for some $ p>1$.

In Section \ref{sec3}, some
 applications of the main results are demonstrated. We   construct confidence intervals for estimating the criticality parameter $\mu$  in terms of $  \frac{Z_{n_0+n}}{Z_{n_0}}  $ and  $ n $; see Propositions \ref{c0kls}  and \ref{zfgsgg}.
The proofs of the main results are given in Sections \ref{sec4} - \ref{sect-last}, by developing a method  different to that used in
   \cite{GLE17} .


\section{Main results} \label{sec2}
\setcounter{equation}{0}

It is well-known that the normalized population size  $$W_n= \frac{Z_n }{\Pi_n},\ \ n \geq 0, $$
 is a nonnegative
martingale  both under the quenched law $\mathbb{P}_\xi$ and under the annealed law $\mathbb{P}$, with respect to
the natural filtration $\mathcal{F}_0=\sigma\{\xi\} ,$
$\mathcal{F}_n=\sigma\{\xi, X_{k,i}, 0\leq k \leq n-1, i\geq 1\}, n\geq 1.$  Then  the limit
 $$W=\lim_{n\rightarrow \infty} W_n$$
 exists $\mathbb{P}$-a.s. by Doob's convergence theorem, and satisfies $\mathbb{E} W \leq 1$ by Fatou's lemma.
 Throughout the paper,   assume that
\begin{equation}\label{dgf103hm}
\mathbb{E} \frac{Z_1}{m_0} \ln^+ Z_1  < \infty.
\end{equation}
Together with the condition that $ p_0(\xi_0)=0 $ a.s.,
condition (\ref{dgf103hm}) implies that
$\mathbb{P}(W >0)=\mathbb{P}(Z_n \rightarrow \infty)=\lim_{n \rightarrow \infty} \mathbb{P}(Z_n >0) =1,$
and that the martingale $W_n$ converges  to $W$ in $\mathbb{L}^1(\mathbb{P})$  (see Athreya and Karlin \cite{AK71b} and also Tanny \cite{T88}).
Clearly, the following decomposition holds:
\begin{equation}\label{decopos}
 \ln Z_n = \sum_{i=1}^n X_i + \ln W_n,
\end{equation}
where $X_i=\ln m_{i-1} (i\geq 1)$ are iid random variables depending only on the environment $\xi.$
The asymptotic behavior of $\ln Z_n$ is crucially affected by  the \emph{associated random walk} 
$$S_n =\sum_{i=1}^n X_i=\ln \Pi_n,  \quad  n\geq 0.$$
By our notation and hypothesis (see \eqref{cond-0}), it follows that
$X=X_1$, $  \mu =\mathbb{E}X >0$ and $ \sigma^2=\mathbb{E}(X-\mu)^2  \in (0, \infty)$; the later implies that
 the random walk $\{ S_n, n\geq 0\}$ is non-degenerate.


 We will need the following Cram\'{e}r  condition on the associated random walk.
\begin{description}
\item[A1.]  The random variable $X=\ln m_0$ has an exponential moment, i.e. there exists a constant $\lambda_0>0$
 such that
\[
\mathbb{E}e^{\lambda_0 X} =\mathbb{E} m_0^{\lambda_0}  < \infty.
\]
\end{description}
Our first result concerns
the  uniform Cram\'{e}r  moderate deviations for 
\begin{equation} \label{def-Zn0n}
Z_{n_0,n}:= \frac{ \ln \frac{Z_{n_0+n}}{Z_{n_0}} - n \mu \ }{  \sigma \sqrt{n}}, \quad n_0\in \mathbb{N}.
\end{equation}

\begin{theorem}\label{th2.2}
Assume condition   \textbf{A1}.  Then  the following results hold  uniformly in $n_0 \in \mathbb N$: for  $n\geq 2$ and   $0 \leq x < \sqrt{  \ln n }, $
\begin{equation}\label{ineq01}
\Bigg| \ln \frac{\mathbb{P}\big(  Z_{n_0,n} \geq x  \big)}{1-\Phi(x)} \Bigg|  \leq C      (1+x^3 ) \frac{   \ln n  }{  \sqrt{n}  };
\end{equation}
for  $n\geq 2$ and   $\sqrt{  \ln n }\leq x = o(\sqrt{n } )  $ as $ n \rightarrow \infty$,
\begin{equation}\label{ineq02}
\Bigg| \ln \frac{\mathbb{P}\big( Z_{n_0,n} \geq x  \big)}{1-\Phi(x)}  \Bigg|
\leq  C   \frac{  x^3   }{  \sqrt{n}  }    .
\end{equation}
The results remain valid when $\frac{\mathbb{P} (  Z_{n_0,n} \geq x  )}{1-\Phi(x)}$ is replaced by  $\frac{\mathbb{P} (  -Z_{n_0,n}\geq x  )}{ \Phi(-x)}$.
\end{theorem}

The uniformity in $n_0$ is interesting in applications. Due to the uniformity,  in \eqref{ineq01}  and (\ref{ineq02})  we can take $n_0$ as  a function of $n$.
Inequality (\ref{ineq02})  coincides with the corresponding result  for  the random walk (cf. \cite{Cramer38} or inequality (1) of \cite{FGL13}), while
in inequality (\ref{ineq01})  there is the additional factor $\ln n$ for BPRE.
When $n_0=0$, the inequalities  (\ref{ineq01}) and  (\ref{ineq02})  but without the factor $\ln n$ have been proved by Grama et al.\ \cite{GLE17} under the additional
condition that  $ \mathbb{E}  \frac{Z_1^{p}}{m_0}   < \infty$ for some $ p>1$.




Theorem  \ref{th2.2}  implies  the following  uniform MDP for $ Z_{n_0,n}$. 
\begin{corollary}\label{corollary02}
Assume  condition  \textbf{A1}.
Let $a_n$ be any sequence of real numbers satisfying $a_n \rightarrow \infty$ and $a_n/\sqrt{n}\rightarrow 0$
as $n\rightarrow \infty$.  Then, for each Borel set $B$,
\begin{eqnarray}
- \inf_{x \in B^o}\frac{x^2}{2} &\leq & \liminf_{n\rightarrow \infty}\frac{1}{a_n^2}\ln \; \inf_{n_0\in \mathbb N} \mathbb{P}\bigg(  \frac{Z_{n_0,n}   }{ a_n  }     \in B \bigg) \nonumber \\
 &\leq& \limsup_{n\rightarrow \infty}\frac{1}{a_n^2}\ln \; \sup_{n_0\in \mathbb N} \mathbb{P}\bigg(\frac{ Z_{n_0,n} }{ a_n  }  \in B \bigg) \leq  - \inf_{x \in \overline{B}}\frac{x^2}{2}   ,   \label{MDP}
\end{eqnarray}
where $B^o$ and $\overline{B}$ denote the interior and the closure of $B$, respectively.
\end{corollary}

%
%
%

The MDP for $Z_{0,n}$ 
has been established by Huang and Liu  \cite{HL12} (see Theorem 1.6 therein)
when the random variable $X= \ln m_0$ satisfies $A_1 \leq m_0$ and $m_0(1+\delta) \leq A^{1+\delta}$ for constants $\delta, A_1$ and $A_2$ satisfying $\delta>0$ and $1< A_1< A$,  and by Wang and Liu  \cite{WL17} under the same condition   \textbf{A1} but in a more general setting.


From  Theorem \ref{th2.2}, using the inequality $|e^y-1| \leq e^C |y|$ valid for $|y| \leq C,$  we obtain the following result about the uniform equivalence to the normal tail.
\begin{corollary}\label{co02s}
Assume condition  \textbf{A1}. Then, uniformly    for $n_0 \in \mathbb N$, as $n\rightarrow \infty$,
\begin{equation} \label{rezero}
\frac{\mathbb{P}\big( Z_{n_0,n} \geq x  \big)}{1-\Phi(x)} =1+ o(1)
\end{equation}
for $x \in [0, \, o(n^{1/6}))$.  The result remains valid when $\frac{\mathbb{P} (  Z_{n_0,n}\geq x  )}{1-\Phi(x)}$ is replaced by  $\frac{\mathbb{P} (  -Z_{n_0,n} \geq x  )}{ \Phi(-x)}$.
\end{corollary}

Inequality (\ref{rezero}) states that the relative error   for normal approximation tends to zero   uniformly for $x \in [0, o(n^{1/6})).$
Notice that the normal range $x \in [0, o(n^{1/6}))$  coincides with the random walk case, under Cram\'{e}r's condition \textbf{A1}. 
In the following Cram\'{e}r moderate deviation result of type Linnik \cite{L61}, we
give a normal range
when the exponential moment condition \textbf{A1} is relaxed to 

\begin{description}
\item[A2.]  The random variable $X=\ln m_0$ has a  sub-exponential moment, i.e.\ there exist  two constants $\lambda_0>0$ and $\gamma \in (0, \frac16]$
 such that
\[
\mathbb{E}\exp\{  \lambda_0  X^{\frac{4 \gamma}{1-2\gamma}} \}<\infty.
\]
\end{description}

\begin{theorem}\label{th03}
Assume condition  \textbf{A2}. Then
\eqref{rezero} holds
uniformly in $n_0\in \mathbb{N}$, for $x \in [0, \, o(n^{\gamma}))$, as $n\rightarrow \infty$. The result  remains valid when $\frac{\mathbb{P} (  Z_{n_0,n}   \geq x  )}{1-\Phi(x)}$ is replaced by   $\frac{\mathbb{P} (  - Z_{n_0,n}  \geq x  )}{ \Phi(- x)}$.
\end{theorem}

Notice that when $\gamma=\frac16,$ Theorem \ref{th03} reduces to  Corollary \ref{co02s}.

We now consider the uniform Berry-Esseen bound for $Z_{n_0,n}$ and $-Z_{n_0,n}$.
The following result under the exponential moment condition  \textbf{A1} can be obtained as a corollary to
Theorem \ref{th2.2}. 

\begin{theorem}\label{cos01}
Assume condition  \textbf{A1}. Then the following holds uniformly in $n_0 \in \mathbb N$:   for $n\geq 2,$
\begin{equation} \label{ths01}
\sup_{x\in \mathbb{R}} \Big|\mathbb{P}\big(  Z_{n_0,n} \leq  x  \big)  -  \Phi(x)   \Big|
\leq  C   \frac{ \ln n   }{  \sqrt{n}  }
\end{equation}
and
\begin{equation}\label{ths02}
\sup_{x\in \mathbb{R}} \Big|\mathbb{P}\big(- Z_{n_0,n} \leq  x  \big)  -  \Phi(x)   \Big|
\leq  C   \frac{ \ln n   }{  \sqrt{n}  }.
\end{equation}
\end{theorem}

In \eqref{ths01} and   \eqref{ths02}  there are the additional factor $\ln n$ for BPRE  compared  to the Berry-Esseen bound  for random walks, for which the best rate is $\frac{ C}{  \sqrt{n}  }$.
We conjecture that the factor $\ln n$ in \eqref{ths01} and   \eqref{ths02} can be removed, just as in the case where $n_0 = 0$ considered in Grama et al.\ \cite{GLE17}.
Actually  Grama et al.\ \cite{GLE17} gave the more general Berry-Esseen bound $ \frac{  C }{   n^{\rho/2}}$    for $Z_{0,n}$ under  a moment condition of order $2+\rho$ on $X$, with $\rho \in (0, 1] $.
We shall prove the same bound for $Z_{n_0,n}$ when $\rho \in (0, \frac{\sqrt{5} -1}{2})$, namely, when the following moment condition holds:
\begin{description}
\item[A3.]  There exists a constant $\rho \in (0, \frac{\sqrt{5} -1}{2})$
 such that
\[
\mathbb{E}X^{2+\rho}    < \infty.
\]
\end{description}
\begin{theorem}\label{th02}
Assume condition  \textbf{A3}. Then uniformly in $n_0 \in \mathbb N$,
\begin{equation}\label{Berr01}
\sup_{x\in \mathbb{R}} \Big|\mathbb{P}\big(  Z_{n_0,n} \leq  x  \big)  -  \Phi(x)   \Big|
\leq  \frac{  C }{   n^{\rho/2}   }
\end{equation}
and
\begin{equation} \label{Berr02}
\sup_{x\in \mathbb{R}} \Big|\mathbb{P}\big(  -Z_{n_0,n}  \leq  x  \big)  -  \Phi(x)   \Big|
\leq   \frac{C   }{   n^{\rho/2}   }.
\end{equation}
\end{theorem}

For $n_0=0$, the inequalities  (\ref{Berr01}) and  (\ref{Berr02})  have been established by Grama et al.\ \cite[Theorem 1.1]{GLE17}  assuming
$\mathbb{E}X^{2+\rho}    < \infty $ for some $\rho \in (0,1]$ and
$ \mathbb{E}  (\frac{Z_1}{m_0})^{p}   < \infty$ for some $p>1$.

\section{Applications to construction of confidence intervals }\label{sec3}
\setcounter{equation}{0}
 Cram\'{e}r moderate deviations  can be applied to constructing  confidence intervals for the  criticality parameter $\mu$.
Assume that $\sigma$ is known, the following two propositions give two  confidence intervals for  $\mu$.
\begin{proposition}\label{c0kls} Assume condition  \textbf{A1}.  Let $\kappa_n \in (0, 1).$  Assume that
\begin{eqnarray}\label{keldet}
 \big| \ln \kappa_n \big| =o \big(  n^{1/3}  \big) .
\end{eqnarray}
Let $$\Delta_n=\frac{\sigma}{ \sqrt{n} }\Phi^{-1}(1-\kappa_n/2)  .$$
Then $[A_n,B_n]$, with
\begin{eqnarray*}
A_n=  \frac{1}{n} \ln \Big(\frac{ Z_{n_0+n} }{ Z_{n_0}} \Big )    -\Delta_n  \quad   \textrm{and} \ \quad  
B_n=  \frac{1}{n} \ln \Big(\frac{ Z_{n_0+n} }{ Z_{n_0}} \Big )  +\Delta_n, 
\end{eqnarray*}
is a $1-\kappa_n$ confidence interval for $\mu$, for $n$ large enough.
\end{proposition}
\emph{Proof.}  By Corollary \ref{co02s},    for $0\leq x=o (n^{1/6}),$
\begin{equation} \label{tphisns4}
\frac{\mathbb{P}(Z_{n_0,n}  \geq x)}{1-\Phi \left( x\right)}=1+o(1)\ \ \  \textrm{and}\ \ \  \frac{\mathbb{P}( Z_{n_0,n}  \leq-x)}{ \Phi \left(- x\right)}=1+o(1).
\end{equation}
 Clearly, the upper $(\kappa_n/2)$th quantile of a standard normal distribution
$$\Phi^{-1}( 1-\kappa_n/2)=-\Phi^{-1}(   \kappa_n/2) = O(\sqrt{| \ln \kappa_n |} ), $$  which, by (\ref{keldet}), is of order $o\big( n^{1/6}  \big).$
Then applying the last equality to (\ref{tphisns4}), we have
\begin{equation}
\mathbb{P}\big(Z_{n_0,n}  \geq  \Phi^{-1}( 1-\kappa_n/2)\big) \sim \kappa_n/2 \ \ \  \textrm{and}\ \ \  \mathbb{P}\big( Z_{n_0,n} \leq -\Phi^{-1}( 1-\kappa_n/2) \big) \sim \kappa_n/2
\end{equation}
as $n\rightarrow \infty.$
Clearly, $Z_{n_0,n}  \leq \Phi^{-1}( 1-\kappa_n/2)$ means that $\mu \geq A_n,   $ while $Z_{n_0,n}  \geq -\Phi^{-1}( 1-\kappa_n/2)$
means  $\mu \leq B_n.  $ This completes the proof of  Proposition \ref{c0kls}. \hfill\qed

When the risk probability $\kappa_n$ goes to $0$, we have the following   result.
\begin{proposition}\label{zfgsgg}  Assume condition  \textbf{A1}. Let $\kappa_n \in (0, 1)$ such that $k_n \rightarrow 0$. Assume that
\begin{eqnarray}\label{3.3sfs}
 \big| \ln \kappa_n \big| =o \big(  n   \big) .
\end{eqnarray}
Let $$\Delta_n=\frac{\sigma}{ \sqrt{n} } \sqrt{ 2 |\ln (\kappa_n/2)|}   .$$
Then $[A_n,B_n]$, with
\begin{eqnarray*}
A_n= \frac{1}{n} \ln \Big(\frac{ Z_{n_0+n} }{ Z_{n_0}} \Big ) -\Delta_n  \quad   \textrm{and} \ \quad  
B_n=  \frac{1}{n} \ln \Big(\frac{ Z_{n_0+n} }{ Z_{n_0}} \Big )  +\Delta_n, 
\end{eqnarray*}
is a $1-\kappa_n$ confidence interval for $\mu$, for $n$ large enough.
\end{proposition}
\emph{Proof.} By Theorem \ref{th2.2}, we have
\begin{equation}\label{ggsg12}
\frac{\mathbb{P}(Z_{n_0,n}  \geq x)}{1-\Phi \left( x\right)}=\exp\bigg\{ \theta_1 C  \frac{(\ln n)^3+x^{3}}{ n^{1/2}} \bigg \}\textrm{\   and \ }\frac{\mathbb{P}(-Z_{n_0,n}  \geq x)}{\Phi \left( -x\right)}=\exp\bigg\{ \theta_2 C  \frac{(\ln n)^3+x^{3}}{ n^{1/2}} \bigg \}
\end{equation}
uniformly for $0\leq x=o (   n^{1/2}  ),$ where $\theta_1, \theta_2 \in [-1, 1]$. Notice that $$1-\Phi \left( x_n\right)=\Phi \left( -x_n\right) \rightarrow \frac{1}{x_n\sqrt{2\pi}}e^{-x_n^2/2}= \exp\bigg\{-\frac{x_n^2}{2}\Big(1+\frac{2}{x_n^2}\ln (x_n\sqrt{2\pi})   \Big) \bigg\} ,\ x_n \rightarrow \infty, $$
and $\gamma \in (0, 1]$. Since $k_n \rightarrow 0$,
the upper $(\kappa_n/2)$th quantile of the distribution $$1-\Big(1-\Phi \left( x\right)\Big)\exp\bigg\{ \theta_1 C  \frac{(\ln n)^3+x^{2+\rho}}{ n^{\rho/2}} \bigg\}$$
converges to  $ \sqrt{ 2 |\ln (\kappa_n/2)|}$, which by (\ref{3.3sfs}) is of order $o\big(n^{1/2}\big)$  as $n\rightarrow \infty.$
Then applying (\ref{ggsg12}) to $Z_{n_0,n} $ and $-Z_{n_0,n}$, by an argument similar to the proof of Proposition \ref{c0kls}, we obtain
the desired result. \hfill\qed


\section{Proof of Theorem \ref{th2.2}} \label{sec4}
\setcounter{equation}{0}
We should prove Theorem \ref{th2.2} for the case of $\frac{\mathbb{P} ( Z_{n_0,n}  \geq x  )}{1-\Phi(x)},$ $ x\geq 0.$ Thanks to the existence a harmonic moment 
(see Lemma  \ref{lemma1}), the case  of $\frac{\mathbb{P} (  -Z_{n_0,n}   \geq x  )}{ \Phi(-x)}$ can be proved in the similar way.
To this end, we start with the proofs of Lemmas   \ref{th200}  and  \ref{th100}, and conclude with the proof of
Theorem \ref{th2.2}. In the sequel, we denote $$\eta_{n,i}= \frac{X_i-\mu}{\sigma\sqrt{n}} , \ \ \ \ \ i=1,...,n_0+n. $$
Then it is easy to see that $\sum_{i=1}^n \mathbb{E} \eta_{n,n_0+i}^2=1.$ Denote $$W_{n_0,n} = \frac{W_{n_0+n}}{W_{n_0}} \ \ \ \textrm{and} \ \ \ W_{n_0,\infty} = \frac{W }{\ \ W_{n_0}}.$$ Then  $(W_{n_0,n})_{n\geq0}$  is also a nonnegative  martingales both under the quenched law $\mathbb{P}_\xi$ and under the annealed law $\mathbb{P}$ with respect to
the natural filtration.

The following lemma gives the upper bound of Theorem \ref{th2.2}.
\begin{lemma} \label{th200}
Assume condition \textbf{A1}. Then the following holds uniformly in $n_0 \in \mathbb N$: \\
 for  $n\geq 2$ and $0 \leq x < \sqrt{  \ln n },$
\begin{equation}\label{ineq3sfd32}
\ln \frac{\mathbb{P}\big(  Z_{n_0,n} \geq x  \big)}{1-\Phi(x)}  \leq C    (1+x^3 )\frac{\ln n }{\sqrt{  n }};
\end{equation}
and for $n\geq 2$ and    $\sqrt{  \ln n }\leq  x = o(\sqrt{n} ),$
\begin{equation}\label{ineq332}
\ln  \frac{\mathbb{P}\big( Z_{n_0,n}  \geq x  \big)}{1-\Phi(x)}  \leq C  \frac{  x^3   }{  \sqrt{n} }    .
\end{equation}
\end{lemma}
\emph{Proof.} We first give a proof for (\ref{ineq332}).
Clearly, by (\ref{decopos}), it holds for  $ \sqrt{ \ln n }  \leq x = o(\sqrt{n} ),$
 \begin{eqnarray}
\mathbb{P}\bigg(  Z_{n_0,n} \geq x  \bigg)
&=&  \mathbb{P}\bigg(   \sum_{i=1}^n \eta_{n,n_0+i} + \frac{\ln W_{n_0,n}}{\sigma\sqrt{n} }\geq x \bigg )
  \leq    \mathbb{P}\bigg(   \sum_{i=1}^n \eta_{n,n_0+i} + \frac{(\ln W_{n_0,n})^+}{\sigma \sqrt{n}} \geq x \bigg ) \nonumber  \\
& \leq&    I_1+I_2,   \label{thn635s}
\end{eqnarray}
where
$$I_1 = \mathbb{P}\bigg(   \sum_{i=1}^n \eta_{n,n_0+i} \geq (x- \frac{  x^{2}}{\sigma \sqrt{n}} )  \bigg ) \ \ \ \textrm{and} \ \ \ I_2= \mathbb{P}\bigg(   \frac{ (\ln W_{n_0,n})^+}{\sigma \sqrt{n} }  \geq \frac{   x^{2}}{\sigma \sqrt{n} } \bigg ).$$
Next, we give some estimations for  $I_1$ and $I_2$.  Notice that $\sum_{i=1}^n \eta_{n,n_0+i}$ is a sum of iid random variables.
 By upper bound of Cram\'{e}r  moderate deviations for sums of iid random variables (cf. inequality (1) of \cite{FGL13}),   we obtain for   $\sqrt{ \ln n } \leq x = o(\sqrt{n} ),$
 \begin{eqnarray*}
I_1 \leq  \bigg(1- \Phi(x-\frac{  x^{2}}{\sigma \sqrt{n} })  \bigg )\exp\bigg\{ \frac{  C}{\sqrt{n} }(x-\frac{ x^{2}}{ \sigma \sqrt{n}})^{3} \bigg\} .
\end{eqnarray*}
Using  the following inequalities
\begin{equation}\label{norb}
\frac 1{\sqrt{2 \pi} ( 1+x)  }e^{- x^2/2}\leq
1-\Phi \left( x\right)   \leq  \frac 1{\sqrt{ \pi} ( 1+x)  }e^{- x^2/2}, \ \   x \geq 0,
\end{equation}
we deduce  that for  $x\geq \ln 2$ and $  \varepsilon_n  \geq 0 $,
\begin{eqnarray}
  \frac{1-\Phi \left( x (1-  \varepsilon_n) \right)}{1-\Phi \left( x\right) } & =& 1+ \frac{ \int_{x (1-  \varepsilon_n)}^x \frac{1}{\sqrt{2\pi}}e^{-t^2/2}dt }{1-\Phi \left( x\right) } \nonumber \\
  &\leq&    1+ \frac{\frac{1}{\sqrt{2\pi}} e^{-x^2(1- \varepsilon_n)^2/2}  x \varepsilon_n  }{ \frac{1}{\sqrt{2 \pi} (1+x)} e^{-x^2/2}  }  \nonumber \\
   &\leq &   1+  C   x^2    \varepsilon_n  e^{   C  x^2 \varepsilon_n }   \nonumber \\
    &\leq&     \exp\Big\{ C    x^2   \varepsilon_n  \Big\}. \label{sfdsh}
\end{eqnarray}
Hence, for   $\sqrt{ \ln n }  \leq x = o(\sqrt{n} ),$
\begin{eqnarray}\label{ines4}
I_1  \leq  \Big(1- \Phi(x )  \Big )\exp\Big\{ C \frac{ x^{3} }{\sqrt{n} }\Big\} .
\end{eqnarray}
By Markov's inequality and (\ref{norb}), it is easy to see that for   $\sqrt{ \ln n } \leq x = o(\sqrt{n} ),$
 \begin{eqnarray}
I_2&=&\mathbb{P}\Big(  W_{n_0,n} \geq \exp\big\{    x^{2}  \big \} \Big)  \nonumber  \\
&\leq & \exp\big\{ -     x^2  \big \}  \mathbb{E} W_{n_0,n} =\exp\big\{ -     x^2  \big \}  \nonumber  \\
&\leq&  C  \frac{ 1+x  }{\sqrt{n} } \Big(1-\Phi(x)\Big)  . \label{ines5}
\end{eqnarray}
Combining (\ref{ines4})  and (\ref{ines5}) together, we obtain for   $\sqrt{ \ln n }\leq  x = o(\sqrt{n} ),$
 \begin{eqnarray*}
\mathbb{P}\bigg( Z_{n_0,n}    \geq x  \bigg )
&\leq& \Big(1- \Phi(x )  \Big )\exp\Big\{ C_1 \frac{ x^{3} }{\sqrt{n} }\Big\} + C_2  \frac{(1+x) }{\sqrt{n} } \Big(1-\Phi(x)\Big)\\
 &\leq& \Big(1- \Phi(x )  \Big )\exp\Big\{ C_3 \frac{ x^{3} }{\sqrt{n} }\Big\},
\end{eqnarray*}
which gives the desired inequality for   $\sqrt{ \ln n }\leq  x = o(n^{1/2}).$

Next, we give a proof for (\ref{ineq3sfd32}).
 When  $1 \leq x < \sqrt{ \ln n }$,  by an argument similar to that of (\ref{thn635s}), we have
  \begin{eqnarray}
\mathbb{P}\big(   Z_{n_0,n}\geq x  \big)
& \leq&    I_3+I_4,  \label{sdfs10}
\end{eqnarray}
 where
 $$I_3 =\mathbb{P}\bigg(   \sum_{i=1}^n \eta_{n,n_0+i} \geq (x- \frac{  x^{2}  \ln n }{\sigma \sqrt{n}} )  \bigg )    \
  \ \ \ \textrm{and} \ \ \ \  I_4= \mathbb{P}\bigg(   \frac{ (\ln W_{n_0,n})^+}{\sigma \sqrt{n} }  \geq \frac{   x^{2} \ln n }{\sigma \sqrt{n} } \bigg ). $$
With   arguments similar to that of (\ref{ines4}) and (\ref{ines5}), we get for   $1 \leq x < \sqrt{ \ln n },$
\begin{eqnarray}
I_3    &\leq&    \bigg(1- \Phi(x-\frac{  x^{2} \ln n }{\sigma \sqrt{n} })  \bigg )\exp\bigg\{ \frac{  C_1}{\sqrt{n} }(x-\frac{ x^{2} \ln n }{ \sigma \sqrt{n}})^{3} \bigg\} \nonumber \\
&\leq& \Big(1- \Phi(x )  \Big )\exp\Big\{ C_2 x^{3} \frac{ \ln n  }{\sqrt{ n} }\Big\}  \label{sdfs11}
\end{eqnarray}
and
\begin{eqnarray}
I_4 &=&\mathbb{P}\bigg(  W_{n_0,n} \geq \exp\Big\{     x^{2} \ln n   \Big \} \bigg)  \nonumber  \\
&\leq & \exp\Big\{ -     x^{2}  \ln n  \Big \}  \mathbb{E} W_{n_0,n} =   \exp\Big\{ -     x^{2}  \ln n  \Big \} \nonumber  \\
& \leq& C \frac{ 1+x}{\sqrt{n} } \Big(1-\Phi(x)\Big)  .  \label{sdfs12}
\end{eqnarray}
Combining (\ref{sdfs10}), (\ref{sdfs11})  and (\ref{sdfs12}) together, we obtain the desired inequality for   $1 \leq x < \sqrt{\ln n } .$

 When  $0 \leq x \leq 1,$  again by an argument similar to that of (\ref{thn635s}), we have
  \begin{eqnarray} \label{sdfs20}
\mathbb{P}\big(  Z_{n_0,n} \geq x  \big)
& \leq&    I_5+I_6,
\end{eqnarray}
 where
 $$I_5 =\mathbb{P}\bigg(   \sum_{i=1}^n \eta_{n,n_0+i} \geq (x- \frac{    \ln n }{\sigma \sqrt{n}} )  \bigg )    \ \ \
  \ \ \ \textrm{and} \ \ \ I_6= \mathbb{P}\bigg(   \frac{ (\ln W_{n_0,n})^+}{\sigma \sqrt{n} }  \geq \frac{   \ln n }{\sigma \sqrt{n} } \bigg ). $$
With   arguments similar to that of (\ref{ines4}) and (\ref{ines5}),
we get for   $0 \leq x \leq 1 ,$
\begin{eqnarray}
I_5    &\leq&    \bigg(1- \Phi(x-\frac{ \ln n }{\sigma \sqrt{n} })  \bigg )\exp\bigg\{ \frac{  C}{\sqrt{n} }(x-\frac{ \ln n }{ \sigma \sqrt{n}})^{3} \bigg\} \nonumber \\
&\leq& \Big(1- \Phi(x )  \Big )\Big(1+ C_1\frac{ \ln n }{\sigma \sqrt{n} } \Big)\Big(1+\frac{ C_2 }{\sigma \sqrt{n} } \Big) \nonumber \\
&\leq& \Big(1- \Phi(x )  \Big ) \Big(1+ C_3 \frac{ \ln n }{  \sqrt{n} } \Big) \label{sdfs21}
\end{eqnarray}
and
\begin{eqnarray}
I_6 &=&\mathbb{P}\bigg(  W_{n_0,n} \geq \exp\Big\{     \ln n   \Big \} \bigg)  \nonumber  \\
&\leq & \exp\Big\{ -      \ln n  \Big \}  \mathbb{E} W_{n_0,n} =  \frac 1n  .  \label{sdfs22}
\end{eqnarray}
Combining  (\ref{sdfs20}), (\ref{sdfs21})  and (\ref{sdfs22}) together, we obtain for   $0 \leq x \leq 1, $
  \begin{eqnarray}
\mathbb{P}\bigg(  Z_{n_0,n} \geq x  \bigg)
& \leq&  \Big(1- \Phi(x )  \Big ) \Big(1+ C_1\frac{ \ln n }{\sigma \sqrt{n} } \Big) + \frac 1n  \nonumber \\
&\leq &  \Big(1- \Phi(x )  \Big ) \Big(1+ C_2 \frac{ \ln n }{  \sqrt{n} } \Big) \nonumber\\
&\leq &  \Big(1- \Phi(x )  \Big ) \exp\Big\{ C_2 \frac{ \ln n }{ \sqrt{n} } \Big\},  \nonumber
\end{eqnarray}
which gives  the desired inequality for   $0\leq x \leq 1$.
This completes the proof of Lemma \ref{th200}.  \hfill\qed

 To prove  the lower bound of Theorem \ref{th2.2}, we shall make use of the following lemma (see Theorem 3.1 of Grama et al.\ \cite{GLE17}).
The  lemma
shows that condition \textbf{A1} implies the existence of harmonic moments of order $a >0.$
\begin{lemma}\label{lemma1}
Assume  condition \textbf{A1}.  There exists a constant $a_0>0$   such that for   $a \in (0, a_0),$
 \begin{eqnarray}
   \mathbb{E}W^{-a}< \infty.
\end{eqnarray}
\end{lemma}

The following lemma gives the lower bound of Theorem \ref{th2.2}.
\begin{lemma}\label{th100}
Assume condition  \textbf{A1}. Then the following holds uniformly in $n_0 \in \mathbb N$: \\
 for  $n\geq 2$ and    $0 \leq x < \sqrt{  \ln n },$
\begin{equation}
\ln \frac{\mathbb{P}\big(  Z_{n_0,n} \geq x  \big)}{1-\Phi(x)} \geq - C  (1+  x^3)  \frac{\ln n }{\sqrt{  n }};
\end{equation}
and for  $n\geq 2$ and  $\sqrt{  \ln n }\leq x = o(\sqrt{n}),$ $n\rightarrow \infty,$
\begin{equation}\label{jnsfg5}
\ln \frac{\mathbb{P}\big( Z_{n_0,n}  \geq x  \big)}{1-\Phi(x)}  \geq- C  \frac{  x^3   }{  \sqrt{n} }    .
\end{equation}
\end{lemma}
\emph{Proof.} We first give a proof for (\ref{jnsfg5}).
Clearly, it holds for   $\sqrt{ \ln n } \leq x = o(\sqrt{n} ),$
 \begin{eqnarray}
\mathbb{P}\bigg(  Z_{n_0,n} \geq x  \bigg)
&=&  \mathbb{P}\bigg(   \sum_{i=1}^n \eta_{n,n_0+i} + \frac{\ln W_{n_0,n}}{\sigma \sqrt{n}} \geq x \bigg )
  \geq     \mathbb{P}\bigg(   \sum_{i=1}^n \eta_{n,n_0+i} \geq x  -\frac{(\ln W_{n_0,n}   )^-}{\sigma \sqrt{n}} \bigg )  \nonumber  \\
& \geq&  \mathbb{P}\bigg(   \sum_{i=1}^n \eta_{n,n_0+i} \geq  x+ \frac{ 4 x^{2}}{a  \sigma \sqrt{n} }   \bigg )  -    \mathbb{P}\bigg(   \frac{ (\ln W_{n_0,n} )^-}{\sigma \sqrt{n}}  \geq \frac{4  x^{2}}{a \sigma \sqrt{n}} \bigg ) \nonumber \\
&=:& P_1-P_2,   \label{ines325}
\end{eqnarray}
where $a$ is a constant satisfying $a \in (0, \min\{a_0,1\})$ with $a_0$ given by Lemma \ref{lemma1}.
Next, we give estimations for terms $P_1$ and $ P_2$.
 By lower bound of Cram\'{e}r moderate deviations for sums of iid random variables (cf. inequality (1) of \cite{FGL13}),  we obtain for   $\sqrt{ \ln n } \leq x = o(\sqrt{n} ),$
 \begin{eqnarray*}
P_1 \geq  \bigg(1- \Phi(x+\frac{ 4 x^{2}}{a \sigma \sqrt{n} })  \bigg )\exp\bigg\{-\frac{  C}{\sqrt{n}}(x+\frac{4x^{2}}{a \sigma \sqrt{n}})^{3} \bigg\} .
\end{eqnarray*}
By an argument similar to that of (\ref{sfdsh}),
we deduce  that for   $x\geq \ln 2$ and $ 0\leq \varepsilon_n \leq 1$,
\begin{eqnarray}\label{th3ds16}
  \frac{1-\Phi \left( x (1+ \varepsilon_n) \right)}{1-\Phi \left( x\right) } \geq     \exp\Big\{ -C   x^2   \varepsilon_n \Big\}.
\end{eqnarray}
Hence, for   $\sqrt{ \ln n } \leq x = o(\sqrt{n} ),$
\begin{eqnarray}\label{inesgd3}
P_1 \geq  \Big(1- \Phi(x )  \Big )\exp\Big\{ -C \frac{ x^{3} }{\sqrt{n} }\Big\} .
\end{eqnarray}
By Markov's inequality,    it is easy to see that for   $\sqrt{ \ln n } \leq x = o(\sqrt{n} ),$
 \begin{eqnarray}
P_2  &=&\mathbb{P}\Big( \ln W_{n_0+n} - \ln W_{n_0} \leq  -    4 x^2 / a   \Big) \nonumber  \\
&=&\mathbb{P}\Big( \ln W_{n_0+n}   \leq  -   2 x^2 / a   \Big)+ \mathbb{P}\Big(  - \ln W_{n_0} \leq  -   2 x^2 / a   \Big) \nonumber \\
&\leq & \exp\Big\{ -  2x^2   \Big \}  \mathbb{E}  W_{n_0+n}  ^{-a  }+ \exp\Big\{ -  2x^2/a    \Big \}  \mathbb{E}  W_{n_0 } \nonumber \\
 &\leq & \exp\Big\{ - 2 x^2   \Big \}  \mathbb{E}  W_{n_0+n}  ^{-a  }+ \exp\Big\{ -  2x^2/a   \Big \}.
\end{eqnarray}
By (\ref{dgf103hm}), it is known that $W_n \rightarrow W$ in $\mathbb{L}^{1}.$ Then we have $W_n=\mathbb{E}[W| \mathcal{F}_n]$ a.s.\ By Jensen's inequality,
we get
 \begin{eqnarray*}
 W_{n_0+n} ^{-a } =(\mathbb{E}[W  | \mathcal{F}_{n_0+n}])^{-a  } \leq \mathbb{E}[  W  ^{-a  } | \mathcal{F}_{n_0+n}].
\end{eqnarray*}
Taking expectations with respect to $\mathbb{P}$ on both sides of the last inequality, we deduce that
  \begin{eqnarray}\label{boudw}
\mathbb{E} W_{n_0+n} ^{-a } \leq \mathbb{E}  W ^{-a  }  .
\end{eqnarray}
By Lemma  \ref{lemma1},  we have for   $\sqrt{ \ln n } \leq x = o(\sqrt{n} ),$
 \begin{eqnarray}
P_2  &  \leq&  \exp\Big\{ - 2 x^2   \Big \}  \mathbb{E}   W^{-a } + \exp\Big\{ -  2x^2    \Big \} \nonumber \\
&\leq& C_1   \exp\Big\{ -  2 x^2  \Big \} \nonumber  \\
& \leq& C_2 \frac{ x}{\sqrt{n}} \Big(1-\Phi(x)\Big) \exp\Big\{ -   x^2  \Big \} . \label{ine455}
\end{eqnarray}
Combining (\ref{ines325}), (\ref{inesgd3}) and (\ref{ine455}) together, we obtain for   $\sqrt{ \ln n } \leq x = o( \sqrt{n}  ),$
 \begin{eqnarray*}
\mathbb{P}\big(   Z_{n_0,n} \geq x  \big)
&\geq& \Big(1- \Phi(x )  \Big )\exp\Big\{ -C_1 \frac{ x^3 }{\sqrt{n}}\Big\} - C_2\frac{ x}{\sqrt{n}} \Big(1-\Phi(x)\Big) \exp\Big\{ -   x^2  \Big \}\\
 &\geq& \Big(1- \Phi(x )  \Big )\exp\Big\{ -C_3 \frac{ x^3 }{\sqrt{n}}\Big\},
\end{eqnarray*}
which gives the desired inequality for     $\sqrt{ \ln n } \leq x = o(\sqrt{n} ).$

For  $0 \leq x < \sqrt{ \ln n },$    the assertion of Lemma \ref{th100}  follows by a similar argument,
but in (\ref{thn635s}) with  $\frac{ 4 x^{2}}{a \sigma \sqrt{n}}   $  replaced by $\frac{ 4 x^{2}  \ln n }{a \sigma \sqrt{n}}$ when $1 \leq x <  \sqrt{ \ln n }$ and $\frac{ 4 x^{2}}{a \sigma \sqrt{n}}   $ replaced by $\frac{ 4  \ln n }{a \sigma \sqrt{n}}$
when  $0 \leq x \leq 1$,   and accordingly in the subsequent statements. Then we get the desired inequality for     $0 \leq x < \sqrt{ \ln n }.$
This completes the proof of Lemma \ref{th100}.
 \hfill\qed

\section{Proof of Corollary \ref{corollary02}}
\setcounter{equation}{0}
We only give a proof for the case of $Z_{n_0,n} . $  The  case of $-Z_{n_0,n} $ can  be  proved in a similar way.
We first show that
\begin{eqnarray}\label{dfgsfdf}
\limsup_{n\rightarrow \infty}\frac{1}{a_n^2}\ln \sup_{n_0\in \mathbb N} \mathbb{P}\bigg(\frac{ Z_{n_0,n} }{ a_n  }  \in B \bigg)   \leq  - \inf_{x \in \overline{B}}\frac{x^2}{2}.
\end{eqnarray}
When $B  =\emptyset,$ the last inequality is obvious, with the convention   $-\inf_{x \in \emptyset}\frac{x^2}{2}=\infty$. Thus,  we may assume that $B  \neq \emptyset.$ Given a Borel set $B\subset \mathbb{R},$ let $x_0=\inf_{x\in B} |x|.$ Clearly, we have $x_0\geq\inf_{x\in \overline{B}} |x|.$
Then, by   Theorem \ref{th2.2},
\begin{eqnarray*}
\sup_{n_0\in \mathbb N} \mathbb{P}\Big(Z_{n_0,n} \in  a_n B \Big)
 &\leq& \sup_{n_0\in \mathbb N} \mathbb{P}\Big(\, \big|Z_{n_0,n}  \big|  \geq a_n x_0\Big)\\
 &\leq&  2\Big( 1-\Phi \left( a_nx_0\right)\Big)
  \exp\bigg\{ C \,  (1+( a_nx_0)^3 ) \frac{ 1+ \mathbf{1}_{ [0,  \sqrt{\ln n})  } (a_nx_0 ) \ln n  }{  \sqrt{n}  }  \bigg\}.
\end{eqnarray*}
Using   (\ref{norb}), after some calculations,
we get
\begin{eqnarray*}
\limsup_{n\rightarrow \infty}\frac{1}{a_n^2}\ln\sup_{n_0\in \mathbb N} \mathbb{P}\bigg(\frac{Z_{n_0,n} \ }{ a_n  }  \in B \bigg)
 \ \leq \  -\frac{x_0^2}{2} \ \leq \  - \inf_{x \in \overline{B}}\frac{x^2}{2} ,
\end{eqnarray*}
which gives (\ref{dfgsfdf}).

Next, we show that
\begin{eqnarray}\label{dfsffsfn}
\liminf_{n\rightarrow \infty}\frac{1}{a_n^2}\ln \inf_{n_0\in \mathbb N} \mathbb{P}\bigg(\frac{Z_{n_0,n} \ }{ a_n  }  \in B \bigg) \geq   - \inf_{x \in B^o}\frac{x^2}{2} .
\end{eqnarray}
When $B^o =\emptyset,$ the last inequality is obvious, with the convention   $ \inf_{x \in  \emptyset}\frac{x^2}{2}=\infty$. Therefore, we may assume that $B^o \neq \emptyset.$
For any given small $\varepsilon_1>0,$ there exists an $x_0 \in B^o,$ such that
\begin{eqnarray*}
 0< \frac{x_0^2}{2} \leq   \inf_{x \in B^o}\frac{x^2}{2} +\varepsilon_1.
\end{eqnarray*}
Since $B^o$ is an open set, for $x_0 \in B^o$ and all small enough $\varepsilon_2 \in (0, |x_0|], $ it holds $(x_0-\varepsilon_2, x_0+\varepsilon_2]  \subset B^o.$
Therefore, $x_0\geq\inf_{x\in  B^o} |x|.$ Without loss of generality, we may assume that $x_0>0.$
Obviously, we have
\begin{eqnarray*}
\inf_{n_0\in \mathbb N} \mathbb{P}\big(Z_{n_0,n}  \in a_n B  \big)   &\geq& \inf_{n_0\in \mathbb N}  \mathbb{P}\big( Z_{n_0,n} \in (a_n ( x_0-\varepsilon_2), a_n( x_0+\varepsilon_2)] \big)\\
&=& \inf_{n_0\in \mathbb N} \Big(   \mathbb{P}\big( Z_{n_0,n} \geq  a_n ( x_0-\varepsilon_2)   \big)-\mathbb{P}\big( Z_{n_0,n}  \geq  a_n( x_0+\varepsilon_2) \big) \Big).
\end{eqnarray*}
Again by Theorem \ref{th2.2}, it is easy to see that
$$\lim_{n\rightarrow \infty} \frac{\sup_{n_0\in \mathbb N} \mathbb{P}\big( Z_{n_0,n} \geq  a_n( x_0+\varepsilon_2) \big) }{\inf_{n_0\in \mathbb N} \mathbb{P}\big( Z_{n_0,n}\geq  a_n ( x_0-\varepsilon_2)   \big)  } =0 .$$
Therefore, by Theorem \ref{th2.2}, it holds for all $n$ large enough,
\begin{eqnarray*}
\inf_{n_0\in \mathbb N} \mathbb{P}\bigg(\frac{ Z_{n_0,n}\ }{ a_n } \in B  \bigg)   &\geq&  \inf_{n_0\in \mathbb N} \frac12 \mathbb{P}\bigg( Z_{n_0,n}    \geq  a_n ( x_0-\varepsilon_2)   \bigg) \\
&\geq& \frac12 \Big( 1-\Phi \left( a_n( x_0-\varepsilon_2)\right)\Big) \\
&&  \times \exp\bigg\{ -C \,  (1+( a_n( x_0-\varepsilon_2))^3 ) \frac{ 1+ \mathbf{1}_{  [0, \sqrt{\ln n})  }(a_n( x_0-\varepsilon_2))  \ln n  }{  \sqrt{n}  }  \bigg\}.
\end{eqnarray*}
Using   (\ref{norb}), after some calculations,
we get
\begin{eqnarray*}
\liminf_{n\rightarrow \infty} \frac{1}{a_n^2}\ln \inf_{n_0\in \mathbb N} \mathbb{P}\bigg(\frac{ Z_{n_0,n} \ }{ a_n  }  \in B \bigg)  \geq  -  \frac{1}{2}( x_0-\varepsilon_2)^2 . \label{ffhms}
\end{eqnarray*}
Letting $\varepsilon_2\rightarrow 0,$  we  deduce that
\begin{eqnarray*}
\liminf_{n\rightarrow \infty}\frac{1}{a_n^2}\ln \inf_{n_0\in \mathbb N}  \mathbb{P}\bigg(\frac{Z_{n_0,n} \ }{ a_n  }  \in B \bigg) \ \geq\ -  \frac{x_0^2}{2}  \  \geq \   -\inf_{x \in B^o}\frac{x^2}{2} -\varepsilon_1.
\end{eqnarray*}
Since   $\varepsilon_1$ can be arbitrarily small, we get (\ref{dfsffsfn}).
Combining (\ref{dfgsfdf}) and (\ref{dfsffsfn}) together, we complete  the proof of Corollary \ref{corollary02}.
 \hfill\qed

\section{Proof of Theorem \ref{th03}}
\setcounter{equation}{0}

 To prove  Theorem \ref{th03}, we shall make use of the following lemma.
\begin{lemma} \label{hkls3}
Assume  condition \textbf{A2}.  There exists a constant $a_0>0$   such that for   $a \in (0, a_0),$
 \begin{eqnarray}
   \mathbb{E}\exp\{a \, |\ln W|^\frac{4 \gamma}{1-2\gamma} \}\mathbf{1}_{\{W\leq 1\}} < \infty.
\end{eqnarray}
\end{lemma}
\emph{Proof.}
Denote
$$   \phi(t) =\mathbb{E} e^{-t W }, $$
for   $t\geq 0$.
From inequality (2.7) of  Grama et al.\ \cite{GLE17}, we have  for   $n\geq1$ and $t \geq KA^n,$
 \begin{eqnarray}\label{fsdsdf}
\phi (t) \leq   \alpha^n + \mathbb{P}(\Pi_n > A^n ),
\end{eqnarray}
where $\alpha \in (0, 1).$  Choose $A$
such that $\ln A > \mu.$ By condition \textbf{A2} and   Theorem 2.1  of \cite{FGL17}, there exists a constant $C>0$ such that  for  all $n \geq 1,$
$$ \mathbb{P}(\Pi_n > A^n) = \mathbb{P}\big( S_n-n\mu > n(\ln A -\mu) \big) \leq \exp\big\{ -C  n^\frac{4 \gamma}{1-2\gamma}\big\}.  $$
From (\ref{fsdsdf}), we get for    $t\geq K A^n,$
 \begin{eqnarray}
\phi (t) \leq  \exp\{ -C  n^\frac{4 \gamma}{1-2\gamma}\}.
\end{eqnarray}
Now for any $t \geq K,$ set $n_0$ be the integer such that $KA^{n_{0}+1} >t\geq KA^{n_{0}}, $ so that
$$    n_0  > \frac{\ln (t/K)}{\ln A} -1   .$$
Then, for any $t \geq K,$
 \begin{eqnarray}
\phi (t) \leq   \exp\{ -C  (\frac{\ln (t/K)}{\ln A} -1 )^\frac{4 \gamma}{1-2\gamma}\} \leq \exp\{ -C_1  (\ln t   )^\frac{4 \gamma}{1-2\gamma}\}.
\end{eqnarray}
By the facts that $\mathbb{P}(W \leq t^{-1}) \leq e \phi (t), t>0,$ and
$$\mathbb{E} \exp\{a \, |\ln W|^\frac{4 \gamma}{1-2\gamma} \}\mathbf{1}_{\{W\leq 1\}}  =\frac{4 a \gamma}{1-2\gamma} \int_1^{\infty} \frac{1}{t}\mathbb{P}(W \leq t^{-1}) (\ln t)^\frac{6 \gamma-1}{1-2\gamma}\exp\{a \, (\ln t)^\frac{4 \gamma}{1-2\gamma} \}dt,$$
it follows that $\mathbb{E} \exp\{a \, |\ln W|^\frac{4 \gamma}{1-2\gamma} \}\mathbf{1}_{\{W\leq 1\}}  < \infty $  for $a \in [0, C_1)$.
\hfill\qed

Now we are in position to prove Theorem \ref{th03}.  We only give a proof of Theorem \ref{th03} for the case of
  $\frac{\mathbb{P} (  Z_{n_0,n} \geq x  )}{1-\Phi(x)}$.  For the case  of $\frac{\mathbb{P} (  -Z_{n_0,n}  \geq x  )}{ \Phi(-x)}$,    Theorem \ref{th03} can be proved in a similar way.
  We first consider the case of $1 \leq x = o( n^\gamma )$. Clearly, it holds for   $1 \leq x = o( n^\gamma ),$
 \begin{eqnarray}
\mathbb{P}\bigg(  Z_{n_0,n} \geq x  \bigg)
&=&  \mathbb{P}\bigg(   \sum_{i=1}^n \eta_{n,n_0+i} + \frac{\ln W_{n_0,n}}{\sigma \sqrt{n}} \geq x \bigg )   \nonumber  \\
& \geq&  \mathbb{P}\bigg(   \sum_{i=1}^n \eta_{n,n_0+i} \geq  x+ \frac{ 2x^{2}}{  \sigma n^{3\gamma }   }   \bigg )  -    \mathbb{P}\bigg(   \frac{  \ln W_{n_0,n} }{\sigma \sqrt{n}}  \geq \frac{2 x^{2}}{  \sigma n^{3\gamma }   } \bigg ) \nonumber \\
&=:& T_1-T_2.   \label{kthn6d35s}
\end{eqnarray}
Next, we give estimations for terms $T_1$ and $ T_2$.
 By lower bound of Linnik type Cram\'{e}r  moderate deviations for sums of iid random variables (cf. Linnik \cite{L61}),  we deduce that  for   $1\leq x = o(n^\gamma ),$
 \begin{eqnarray*}
T_1 \geq  \bigg(1- \Phi(x+\frac{2 x^{2}}{  \sigma n^{3\gamma }   })  \bigg )\Big(1-o(1)\Big).
\end{eqnarray*}
Hence, by (\ref{th3ds16}), we get for   $1\leq x = o(n^\gamma ),$
\begin{eqnarray}
T_1 &\geq&  \Big(1- \Phi(x )  \Big )\exp\Big\{ -C \frac{ x^{3} }{ n^{3\gamma} }\Big\}\Big(1-o(1)\Big) \nonumber \\
&=&\Big(1- \Phi(x )  \Big )\Big(1-o(1)\Big) . \label{kinesgd3}
\end{eqnarray}
By Markov's inequality,    it is easy to see that for   $1\leq x = o(n^\gamma ),$
 \begin{eqnarray*}
T_2   &=&\mathbb{P}\Big(  \ln W_{n_0+n}   - \ln W_{n_0}   \geq    2  x^2 n^{\frac12-3\gamma}   \Big)  \nonumber \\
 &\leq&\mathbb{P}\Big(  \ln W_{n_0+n}      \geq    x^2 n^{\frac12-3\gamma}   \Big)+\mathbb{P}\Big(    - \ln W_{n_0}   \geq      x^2 n^{\frac12-3\gamma}   \Big)  \nonumber \\
  &\leq&  \exp\Big\{ - \frac{a_0}2 ( x^2 n^{\frac12-3\gamma})^{\frac{4 \gamma}{1-2\gamma}}   \Big \}  \mathbb{E}  \exp\{  \frac{a_0}2 \, |\ln W_{n_0+n}|^{\frac{4 \gamma}{1-2\gamma}} \}\mathbf{1}_{\{W_{n_0+n} \leq 1\}} \nonumber \\
  && + \exp\Big\{ - \frac{a_0}2 ( x^2 n^{\frac12-3\gamma})^{\frac{4 \gamma}{1-2\gamma}}   \Big \}  \mathbb{E}  \exp\{  \frac{a_0}2 \, |\ln W_{n_0 }|^{\frac{4 \gamma}{1-2\gamma}} \}\mathbf{1}_{\{W_{n_0 } \leq 1\}} .
\end{eqnarray*}
Recall that $W_n=\mathbb{E}[W| \mathcal{F}_n]$ a.s. Since $f(x)= \exp\{  \frac{a_0}2 \, |\ln x|^{\frac{4 \gamma}{1-2\gamma}} \}\mathbf{1}_{\{x\leq 1\}}$ is
 convex in $(0, 1],$ by Jensen's inequality,
we get
 \begin{eqnarray*}
  f(W_n) =f(\mathbb{E}[W| \mathcal{F}_n])  \leq\mathbb{E}[f(W)| \mathcal{F}_n].
\end{eqnarray*}
Taking expectations with respect to $\mathbb{P}$ on both sides of the last inequality, we deduce that
  \begin{eqnarray*}
\mathbb{E}[\exp\{  \frac{a_0}2 \, |\ln W_n|^{\frac{4 \gamma}{1-2\gamma}} \} \mathbf{1}_{\{W_n\leq 1\}}]  \leq \mathbb{E}[\exp\{  \frac{a_0}2 \, |\ln W|^{\frac{4 \gamma}{1-2\gamma}} \}\mathbf{1}_{\{W\leq 1\}} ]  .
\end{eqnarray*}
By Lemma  \ref{lemma1},  we have for   $1\leq x = o(n^\gamma ),$
 \begin{eqnarray}
T_2  &  \leq&  C \exp\Big\{ - \frac{a_0}2 ( x^2 n^{\frac12-3\gamma})^\frac{4 \gamma}{1-2\gamma}  \Big \}  \nonumber \\
&\leq&   \frac{C_1 }{\sqrt{n}} \Big(1-\Phi(x)\Big)   .   \label{kine455}
\end{eqnarray}
Combining (\ref{kthn6d35s}), (\ref{kinesgd3}) and (\ref{kine455}) together, we obtain for   $1\leq x = o(n^\gamma ),$
 \begin{eqnarray}
\mathbb{P}\big(  Z_{n_0,n} \geq x  \big)
&\geq& \big(1- \Phi(x )  \big )\big(1-o(1)\big) -\frac{C}{\sqrt{n}} \big(1-\Phi(x)\big) \nonumber \\
 &\geq& \big(1- \Phi(x )  \big )\big(1-o(1)\big).  \label{tnd012}
\end{eqnarray}
 Similarly, we can prove that for   $1\leq x = o(n^\gamma ),$
 \begin{eqnarray}
\mathbb{P}\big(  Z_{n_0,n} \geq x  \big)
 &\leq& \big(1- \Phi(x )  \big )\big(1+o(1)\big). \label{tnd013}
\end{eqnarray}
Combining (\ref{tnd012}) and (\ref{tnd013}) together,  we have  for   $1\leq x = o(n^\gamma ),$
 \begin{eqnarray*}
\mathbb{P}\big(  Z_{n_0,n} \geq x  \big)
 &=& \big(1- \Phi(x )  \big )\big(1+o(1)\big).
\end{eqnarray*}
This completes the proof of Theorem \ref{th03}.
  For $0\leq x \leq 1,$ Theorem \ref{th03} can be proved in a similar way, but in (\ref{kthn6d35s}) with  $\frac{ 2x^{2}}{  \sigma n^{3\gamma }   } $  replaced by $\frac{ 2 }{  \sigma n^{3\gamma }   } $,  and accordingly in the subsequent statements.
 \hfill\qed

\section{Proof of Theorem \ref{cos01}}
\setcounter{equation}{0}
We only give a proof of (\ref{ths01}). Inequality (\ref{ths02}) can  be  proved in a similar way.
Clearly, it holds
\begin{eqnarray}
&& \sup_{ x \in \mathbb{R}   }  \big|\mathbb{P}\big(   Z_{n_0,n} \leq x  \big) -  \Phi \left( x\right) \big| \nonumber  \\
 &&\ \ \ \ \ \leq  \sup_{   x > n^{1/8}} \big|\mathbb{P}\big(  Z_{n_0,n} \leq x  \big) -  \Phi \left( x\right) \big|    + \sup_{  0 \leq x \leq  n^{1/8} } \big|\mathbb{P}\big(  Z_{n_0,n} \leq x  \big) -  \Phi \left( x\right) \big| \nonumber\\
&& \ \ \ \ \ +  \sup_{ - n^{1/8}\leq x \leq 0 } \big|\mathbb{P}\big(  Z_{n_0,n} \leq x  \big) -  \Phi \left( x\right) \big|  +\sup_{   x < - n^{1/8}} \big|\mathbb{P}\big(  Z_{n_0,n} \leq x  \big) -  \Phi \left( x\right) \big| \nonumber  \\
&&\ \ \ \ \ =:  H_1 + H_2+H_3+H_4.   \label{ineq0d10}
\end{eqnarray}
By Theorem \ref{th2.2} and (\ref{norb}), it is easy to see that
\begin{eqnarray*}
H_1 &= & \sup_{   x > n^{1/8}} \big|\mathbb{P}\big(  Z_{n_0,n} > x  \big)- \big(1 -  \Phi \left( x\right) \big) \big| \\
 &\leq & \sup_{   x > n^{1/8}}  \mathbb{P}\big(  Z_{n_0,n} > x  \big) + \sup_{   x > n^{1/8}}  \big(1 -  \Phi \left( x\right) \big)   \\
 &\leq &   \mathbb{P}\big(  Z_{n_0,n} > n^{1/8}   \big) +    \big(1 -  \Phi \left( n^{1/8} \right) \big) \\
 &\leq &  \big(1 -  \Phi \left( n^{1/8} \right) \big)e^C   +    \exp\{ -\frac12 n^{1/4}\} \\
 & \leq &  C_1 \frac{\ln n}{\sqrt{n}}
\end{eqnarray*}
and
\begin{eqnarray*}
H_4   &\leq & \sup_{   x <- n^{1/8}}  \mathbb{P}\big(  Z_{n_0,n} \leq x  \big) + \sup_{   x <- n^{1/8}}     \Phi \left( x\right)     \\
 &\leq &   \mathbb{P}\big(  Z_{n_0,n} \leq -n^{1/8}   \big) +     \Phi \left( -n^{1/8} \right)  \\
 &\leq &   \Phi \left(- n^{1/8} \right)  e^C   +    \exp\{ -\frac12 n^{1/4}\} \\
 & \leq &  C_2 \frac{\ln n}{\sqrt{n}}.
\end{eqnarray*}
By Theorem \ref{th2.2} and the inequality $|e^x-1|\leq |x|e^{|x|},$ we have
\begin{eqnarray*}
H_2 &=&\sup_{0\leq   x \leq n^{1/8}} \big|\mathbb{P}\big(  Z_{n_0,n} > x  \big)- \big(1 -  \Phi \left( x\right) \big) \big| \nonumber \\
  &\leq & \sup_{0\leq   x \leq n^{1/8}} \big(1-\Phi(x) \big) \big| e^{   C  ( 1+ x^3) (\ln n)/ \sqrt{n} }   -1 \big| \nonumber\\
 &\leq & C_3 \frac{\ln n}{\sqrt{n}}
 \end{eqnarray*}
 and
 \begin{eqnarray*}
H_3 &=& \sup_{ - n^{1/8}\leq x \leq 0 } \big|\mathbb{P}\big(  Z_{n_0,n} \leq x  \big) -  \Phi \left( x\right) \big| \nonumber \\
  &\leq & \sup_{ - n^{1/8}\leq x \leq 0 } \Phi(x) \big| e^{   C  ( 1+ |x|^3) (\ln n)/ \sqrt{n} }   -1 \big| \nonumber\\
 &\leq & C_4 \frac{\ln n}{\sqrt{n}}.
 \end{eqnarray*}
Applying  the bounds  of $H_1, H_2, H_3$ and $H_4$ to (\ref{ineq0d10}),   we obtain  inequality (\ref{ths01}). This completes the proof of Theorem \ref{cos01}. \hfill\qed

\section{Proof of Theorem \ref{th02}} \label{sect-last}
\setcounter{equation}{0}
We should prove Theorem \ref{th02} for the case of $ Z_{n_0,n}. $
The cases of $ -Z_{n_0,n}$ can be proved in the similar way.
To prove  the lower bound of Theorem \ref{th02}, we shall  make use of the following lemma, which
is an improvement on Lemma 2.3 of Grama et al.\ \cite{GLE17}, in which   $p \in (0, 1+  \rho/2)$ instead of $p \in (0, 1+  \rho)$.

\begin{lemma}\label{lmfdgnb}
Assume condition \textbf{A3}. Then  for   $p \in (0, 1+  \rho),$
\begin{equation}
  \mathbb{E}|\ln W |^p< \infty \ \ \  \textrm{and}\ \ \ \sup_{n \in \mathbb{N}} \mathbb{E}|\ln W_n|^p< \infty.
\end{equation}
\end{lemma}
\emph{Proof.} By Jensen's inequality, it is enough to prove   Lemma \ref{lmfdgnb} for $p \in [1, 2+\rho).$
From (2.7) of  of Grama et al.\ \cite{GLE17}, we have  for all $n\geq1$ and $t \geq KA^n,$
 \begin{eqnarray}\label{fssdf}
\phi (t) \leq   \alpha^n + \mathbb{P}(\Pi_n > A^n ),
\end{eqnarray}
where $\alpha \in (0, 1).$ Recall that  $\mu=\mathbb{E}X$ and $S_n=\ln \Pi_n = \sum_{i=1}^nX_i.$ Then $ S_n$ is a sum of  iid random variables with $(2+\rho)$-moments. Choose $A$
such that $\ln A > \mu.$ By Nagaev's inequality (see Corollary 1.8 of Nagaev \cite{N79} or Corollary 2.5 of \cite{FGL17}), there exists a constant $C>0$ such that  for $n \in \mathbb{N},$
$$ \mathbb{P}(\Pi_n > A^n) = \mathbb{P}\big( S_n-n\mu > n(\ln A -\mu) \big) \leq  \frac{C}{n^{1+\rho}} .  $$
From (\ref{fssdf}), we get for all $n$ large enough and $t\geq K A^n,$
 \begin{eqnarray}
\phi (t) \leq  \frac{C}{n^{1+\rho}}.
\end{eqnarray}
Now for any $t \geq K,$ set $n_0$ be the integer such that $KA^{n_{0}+1} >t\geq KA^{n_{0}}, $ so that
$$    n_0  > \frac{\ln (t/K)}{\ln A} -1   .$$
Then, for any $t \geq K,$
 \begin{eqnarray}
\phi (t) \leq   C_0 ( \ln t)^{-1-\rho}.
\end{eqnarray}
By  the facts that $\mathbb{P}(W \leq t^{-1}) \leq e \phi (t), t>0,$ and
$$\mathbb{E} |\ln  W|^p \mathbf{1}_{\{  W \leq 1 \} }  =p \int_1^{\infty} \frac{1}{t}(\ln t)^{p-1}\mathbb{P}(W \leq t^{-1}) dt,  $$
it follows that $\mathbb{E}|\ln W |^p\mathbf{1}_{ \{W \leq 1 \}}   < \infty $  for $p \in [1, 1+\rho)$.
Using the inequality $|\ln x |^p  \leq C x, x> 1$, we deduce that  $\mathbb{E}|\ln W |^p\mathbf{1}_{\{W > 1\}}  \leq C \mathbb{E} W   \leq C \mathbb{E} W_n=C.$
Thus, we have $$\mathbb{E}|\ln W |^p = \mathbb{E}|\ln W |^p\mathbf{1}_{ \{W \leq 1 \}} + \mathbb{E}|\ln W |^p\mathbf{1}_{\{W >1\}}  < \infty.  $$ Notice that $x \mapsto |\ln x|^p\mathbf{1}_{\{0<x \leq 0\}}$ is a non-negative and convex function for $p \in [1, 1+\rho).$ By Lemma 2.1 of Huang and Liu  \cite{HL12},
we have
$ \sup_{n}   \mathbb{E}|\ln W_n |^p \mathbf{1}_{\{W_n \leq 1\}}  =   \mathbb{E}|\ln W  |^p \mathbf{1}_{\{W \leq 1\}}  < \infty.  $
It is also easy to see that for $p \in [1, 1+\rho),$
\begin{eqnarray*}
  \sup_{n}   \mathbb{E}|\ln W_n |^p &=&  \sup_{n}   \mathbb{E}|\ln W_n |^p \mathbf{1}_{\{W_n \leq 1\}} +  \sup_{n}   \mathbb{E}|\ln W_n |^p \mathbf{1}_{\{W_n > 1\}}  \\
   &\leq&   \mathbb{E}|\ln W  |^p \mathbf{1}_{\{W \leq 1\}} +C \mathbb{E}  W_n  = \mathbb{E}|\ln W  |^p \mathbf{1}_{\{W \leq 1\}} +C  < \infty.
\end{eqnarray*}
This completes the proof of Lemma \ref{lmfdgnb}.
\hfill\qed

Now we are in position to prove Theorem \ref{th02}.
We first prove that for   $x\in \mathbb{R},$
\begin{equation}\label{fs525}
\mathbb{P}\Big(  Z_{n_0,n} \leq x  \Big)  -  \Phi(x)    \leq \frac{C}{n^{\delta/2}}.
\end{equation}
It is easy to see that
 \begin{eqnarray}
\mathbb{P}\bigg(  Z_{n_0,n} \leq x  \bigg)
  \leq  \mathbb{P}\bigg(   \sum_{i=1}^n \eta_{n,n_0+i} -\frac{ (\ln W_{n_0,n})^- }{\sigma \sqrt{n}} \leq x \bigg )
  \leq     R_1+R_2,  \label{Berryup}
\end{eqnarray}
where
$$R_1 = \mathbb{P}\bigg(   \sum_{i=1}^n \eta_{n,n_0+i} \leq  x + \frac{ 2}{ \sigma n^{\rho/2} }   \bigg ) \ \ \ \textrm{and} \ \ \ R_2= \mathbb{P}\bigg(    \frac{ (\ln W_{n_0,n})^-}{\sigma \sqrt{n} }  \geq \frac{ 2}{\sigma n^{\rho/2} }  \bigg ).$$
Next, we give estimations for   $R_1$ and $R_2$. By the Berry-Esseen bound for a sum of iid random variables,  we obtain
 \begin{eqnarray}
R_1 &\leq&  \Phi(x+ \frac{ 2}{\sigma n^{\rho/2} }) + \frac{C_1}{n^{\rho/2}} \nonumber \\
    &\leq&  \Phi(x) + \frac{C_2 }{n^{\rho/2}}   .
\end{eqnarray}
Notice that when $\rho \in (0, (\sqrt{5} -1)/2  ),$  we have $p:=\frac{\rho}{1-\rho} < 1+  \rho$.
By Markov's inequality and Lemma \ref{lmfdgnb}, it is easy to see that
 \begin{eqnarray}
R_2&=&\mathbb{P}\Big( |\ln( W_{n_0+n} / W_{n_0}) |   \geq   2 n^{(1-\rho)/2}   \Big)  \nonumber  \\
&\leq&\mathbb{P}\Big( |\ln  W_{n_0+n}|+| \ln W_{n_0}  |   \geq   2 n^{(1-\rho)/2}   \Big)  \nonumber  \\
&\leq&\mathbb{P}\Big( |\ln  W_{n_0+n}|   \geq     n^{(1-\rho)/2}   \Big) + \mathbb{P}\Big(  | \ln W_{n_0}  |   \geq    n^{(1-\rho)/2}   \Big)  \nonumber  \\
&\leq & n^{-p(1-\rho)/2} \mathbb{E} |\ln  W_{n_0+n}|^p +   n^{-p(1-\rho)/2} \mathbb{E} |\ln  W_{n_0}|^p \leq 2n^{-\rho/2} \sup_n \mathbb{E} |\ln  W_n|^p   \nonumber  \\
&\leq&   \frac{C }{n^{\rho/2}}  .
\end{eqnarray}
Applying the upper bounds of $R_1$ and $R_2$  to  (\ref{Berryup}), we obtain  (\ref{fs525}).

Next, we prove that for   $x\in \mathbb{R},$
\begin{equation}\label{fs525d}
\mathbb{P}\big(  Z_{n_0,n} \leq x  \big)  -  \Phi(x)  \geq - \frac{C}{n^{\delta/2}}.
\end{equation}
Clearly, it holds
 \begin{eqnarray}
\mathbb{P}\bigg(  Z_{n_0,n} \leq x  \bigg)
  \geq  \mathbb{P}\bigg(   \sum_{i=1}^n \eta_{n,n_0+i} +\frac{ (\ln W_{n_0,n})^+ }{\sigma \sqrt{n}} \leq x \bigg )
  \geq     R_3-R_4,  \label{Berrylow}
\end{eqnarray}
where
$$R_3 = \mathbb{P}\bigg(   \sum_{i=1}^n \eta_{n,n_0+i} \leq  x - \frac{ 1}{ \sigma n^{\rho/2} }   \bigg ) \ \ \ \textrm{and} \ \ \ R_4= \mathbb{P}\bigg(    \frac{ (\ln W_{n_0,n})^+}{\sigma \sqrt{n} }  \geq \frac{ 1}{\sigma n^{\rho/2} }  \bigg ).$$
Again by the Berry-Esseen bound for a sum of iid random variables,  we obtain
 \begin{eqnarray}
R_3 &\geq&  \Phi(x-\frac{ 1}{\sigma n^{\rho/2} }) - \frac{C_1}{n^{\rho/2}} \nonumber \\
    &\geq&  \Phi(x) - \frac{C_2 }{n^{\rho/2}}   .
\end{eqnarray}
Again by Markov's inequality, we get
 \begin{eqnarray}
R_4&\leq&\mathbb{P}\Big(  W_{n_0,n}   \geq  \exp\{   n^{(1-\rho)/2} \}  \Big)  \nonumber  \\
&\leq &   \exp\{ -  n^{(1-\rho)/2} \}  \mathbb{E}  W_{n_0,n}     \leq  \exp\{ -  n^{(1-\rho)/2} \}  \nonumber  \\
&\leq&   \frac{C }{n^{\rho/2}}  .
\end{eqnarray}
Applying the upper bounds of $R_3$ and $R_4$  to  (\ref{Berrylow}), we obtain  (\ref{fs525d}).

Combining (\ref{fs525}) and (\ref{fs525d}) together, we get
\begin{equation}
\Big| \mathbb{P}\big(  Z_{n_0,n} \leq x  \big)  -  \Phi(x)  \Big|  \leq \frac{C}{n^{\delta/2}},
\end{equation}
which gives  the desired inequality.    \hfill\qed


\end{document}